\let\documentclass\relax\documentclass % Для тупого арХива
\input AHTOHFIE.STY
\def\coauthor{Misha}

\UDC{\hskip-0.1em
519.157.1%   %Задачи о покрытии, о минимальных системах представителей
+519.175.1%   %Изоморфизм графов. Симметрии графов
+519.176%     %Экстремальные задачи теории графов
%+512.544.42% %Группы автоморфизмов и представления групп в группах
%             %автоморфизмов алгебраических систем
%+512.544.44% %Характеристические и вполне характеристические подгруппы
%+512.544.7%   %Аппроксимация групп
}

\title{
Invariant covers of multipartite hypergraphs
}
\author{%
Anton A. Klyachko%$^{\flat\sharp}$
\quad
Mikhail S. Terekhov%$^\flat$
}
\address{
%$^\flat$%
\myAddressWC
\qquad\qquad
kombox.ver.2.0@yandex.ru
}

\grants{\RSF 22-11-00075}

\abstract{%
We prove the following ``symmetric analogue" of Lov\'asz's
estimate (1975): if an $r$-partite hypergraph of rank $r\ge2$ has a cover
of cardinality $n<\infty$, then it admits a cover of cardinality at
most $nr/2$, which is invariant with respect to all automorphisms
preserving the parts. We obtain also symmetric analogues of
generalisations of Lov\'asz's estimate due to Aharoni, Holzman, and
Krivelevich (1996).
}

%%%%%%%%%%%%%%
\s 1.
Introduction

A \emph{hypergraph} is a pair $\Gamma=(V,E)$, where $V=V(\Gamma)$ is
a set (of \emph{vertices}), and
$E=E(\Gamma)$ is a set of subsets of $V$;
the elements of $E=E(\Gamma)$ are called \emph{edges} (or
\emph{hyperedges}).  The \emph{rank} of a hypergraph $\Gamma$ is the
cardinal~$\rk\Gamma=\sup\limits_{e\in E(\Gamma)}|e|$ (we are
interested only in finite-rank hypergraphs).
A set $X\subseteq V(\Gamma)$ is a \emph{\(vertex\) cover}
if it intersects each edge. The minimal cardinality of a
vertex cover of $\Gamma$ is denoted by
$\tau(\Gamma)$.
We call $\Gamma$ \emph{finitely covered} if
$\tau(\Gamma)<\infty$.
%(или, что то же самое, $\tau^*(\Gamma)<\infty$).

If a group $G$ acts by automorphisms on a hypergraph $\Gamma$,
then $\tau_G(\Gamma)$ denotes the minimal cardinality of a $G$-invariant
cover.
For instance, take a cube
$V(\Gamma)=\{0,1\}^3$ and
\-
$E(\Gamma)=\{\hbox{6 faces of the cube}\}=
\bigl\{\{(x_1,x_2,x_3)\in V(\Gamma) \mid x_i=\epsilon\}
\;\bigm|\; i=1,2,3;\ \epsilon=0,1\bigr\}$,
then $\tau(\Gamma)=2$ (the set $\{(0,0,0),\;(1,1,1)\}$ is a minimal
cover), while $\tau_{\Aut\Gamma}=8$ (as the hypergraph
is vertex-transitive);
\-
if we take $\{\hbox{12 edges of the cube}\}$ as $E(\Gamma)$,
then $\tau(\Gamma)$ become four, while
$\tau_{\Aut\Gamma}$ remains eight.

\enditem
It is known [KL21] that, for finitely covered
finite-rank hypergraphs,
$\tau_G(\Gamma)\le\tau(\Gamma)\cdot\rk\Gamma$.
In the general case, this estimate is sharp, but in some
interesting special cases, the ``cost of symmetry" turned out to be
lower, see [T22], [KT25], and [T25]. For example, the
Dulmage--Mendelsohn decomposition
[DM58] almost immediately implies that, if~$\Gamma$ is a
biparite graph, and the action of $G$ preserves parts, then
$\tau_G(\Gamma)=\tau(\Gamma)$, i.e., the symmetry in this case is freely
available. Our purpose is to extend this observation
(playing a substantial role in [KT25]) to multipartite hypergraphs.

The main idea of this paper is that
\dispno{\it
relations between covers and invariant covers are
similar to relations between fractional covers and covers.
}(*)%
Recall that a \emph{fractional \(vertex\) cover} of a hypergraph $\Gamma$
is a function $t$ from $V(\Gamma)$ to the unit
interval $[0;1]\subset\R$ such that $\sum\limits_{v\in e}t(v)\ge1$ for each
edge $e$. The \emph{size} of a fractional cover $t$ is
$|t|\:=\sum\limits_{v\in V(\Gamma)}t(v)$ (where the infinite sum means
the supremum of all finite subsums). The infimum of sizes of fractional
covers of~$\Gamma$ is denoted by $\tau^*(\Gamma)$. Certainly,
$\tau^*(\Gamma)\le\tau(\Gamma)$ (because covers
correspond to fractional covers taking value in~$\{0,1\}$).
There are many papers on relations between
$\tau^*(\Gamma)$ and $\tau(\Gamma)$ for finite hypergraphs:
see, e.~g., [L75], [F\"u88], and [AHK96]. For
instance,
\dispno{\sl
$
\tau(\Gamma)/\tau^*(\Gamma)\le r/2
$
for any finite $r$-partite hypergraph of rank $r\ge2$.}
(\hbox{\emph{Lov\'asz's estimate} [L75]})

\vglue-10pt

Covers of infinite hypergraphs has not been so extensively studied, but
see, e.g., [A85], [AH92], [BZ19], and [Ko19].
We do not assume
the finiteness of hypergraphs in this paper, but it turns out that
the problem under consideration is reduced to the finite case,
see the theorem on finite trace in Section~2. In particular,
Lov\'asz's estimate and its generalisations due to Aharoni, Holzman,
and Krivelevich [AHK96] (see Section 2) turn out to be valid for
infinite hypergraphs.

The main result of this paper is a symmetric analogue of
Aharoni--Holzman--Krivelevich's bound (see Section 3).
The following is the simplest special
case of this result.

\proclaim Symmetric analogue of Lov\'asz's estimate.
If, on an $r$-partite finitely covered hypergraph $\Gamma$ of rank $r\ge2$,
%\newline
a group~$G$ acts
%\newline
preserving the parts, then
$
\tau_G(\Gamma)/\tau(\Gamma)\le r/2
$
\(and the inequality is strict for $r\ge3$\).

This serves as an illustration of thesis $(*)$, but there is a nuance:
Lov\'asz's estimate is sharp [AHK96], while the sharpness of the symmetric
analogue is an open quetion, see the last section.

Another illustration of $(*)$ (as we understand now) is 
the bound from [KL21] mentioned above:
$\tau_G(\Gamma)/\tau(\Gamma)\le\rk\Gamma$
for finitely covered finite-rank hypergraphs,
which is a ``symmetric analogue'' of the estimate
$\tau(\Gamma)/\tau^*(\Gamma)<\rk\Gamma$
for $\rk\Gamma\ge2$ [L76].
There is again a nuance here: the bound from [L76] is strict,
while the inequality from [KL21] becomes an equality very often
(see [T22], [KT25], and [T25]).

%%%%%%%%%%%%%%%%%%%%%%%%%%%%%%%%
\s
2. Finite approximations

It is easy figure out that
$\tau^*(\Gamma)$ and $\tau(\Gamma)$ are finite or infinite simultaneously
([AH92], Lemma 1.6).
If all edges of a hypergraph $\Gamma$ are finite, then
the infimum $\tau^*(\Gamma)$ is attained ([AH92], Theorem 2.1),
i.e., there exists a fractional cover of size $\tau^*(\Gamma)$.
For finite-rank hypergraphs, a stronger fact holds.

\proclaim Finite-fractional-cover lemma.
For each finite-size fractional cover $t$ of a finite-rank hypergraph,
some finite cut~$t[K]$ of $t$ is a fractional cover too.  In particular,
the infimum $\tau^*(\Gamma)$ is attained only on fractional covers with
finite supports.
\newline\rm
Henceforth, a \emph{finite cut} $t[K]$ of a fractional cover $t$ is the
function $V(\Gamma)\to[0;1]$ coinciding with $t$ on a finite
set~$K\subseteq V(\Gamma)$ and vanishing on $V\setminus K$.

\Proof
Induction on $r=\rk(\Gamma)$.
The set $U=\{u\in V(\Gamma)\mid t(u)\ge1/r\}$ is finite
(because the size of $t$ is finite) and is a cover
(because $t$ is a fractional cover).
For each $u\in U$ such that $t(u)\ne1$, consider the
hypergraph~$\Gamma_u$ with vertices
$V(\Gamma_u)=V(\Gamma)\setminus\{u\}$
and edges
$E(\Gamma_u)=\{e\setminus\{u\} \mid u\in e\in E(\Gamma)\}$.
\-
The function ${1\over1-t(u)}t$ is a fractional cover of this graph.
Indeed, for each $e_u\in E(\Gamma_u)$,
we have $e_u=e\setminus\{u\}$, where $e\in E(\Gamma)$, and
$
\sum\limits_{v\in e_u}{1\over1-t(u)}t(v)
=
{1\over1-t(u)}\cdot
\Bigl(
-t(u)+
\underbrace{\sum\limits_{v\in e}t(v)}_
{\ge1\rlap{\hbox{\small(because $\scriptstyle t$ is a fractional cover)}}}
\Bigr)
\ge
{1\over1-t(u)}\cdot\(-t(u)+1\)
\ge1
$;
\-
and $\rk\Gamma_u<\rk\Gamma$.

\enditem
By the induction hypothesis, a
finite cut ${1\over1-t(u)}t[K_u]$ is a fractional cover
of $\Gamma_u$.
Therefore,
$$
v\mapsto\max\{t[K_u\cup L](v)\mid u\in U,\; t(u)\ne1\},
\qbox{where $L=\{u\in V(\Gamma)\mid t(u)=1\}$,}
$$
is a finite cut of $t$ and a
fractional cover of $\Gamma$. This completes the proof.

\smallskip

This lemma makes it possible to extend results on finite hypergraphs to
infinite ones as follows.
A \emph{finite trace} of a hypergraph $\Gamma$ is a hypergraph
$\Gamma[K]$, whose vertices are a finite set $K$ of vertices
of $\Gamma$, and
$E(\Gamma[K])=\{e\cap K \mid e\in E(\Gamma)\}$.

\proclaim Finite-trace theorem.
For any
%04.02.26
finitely covered
finite-rank hypergraph $\Gamma$,
there exists a finite trace $\Gamma[K]$ such that
$\tau(\Gamma)=\tau(\Gamma[K])$ and $\tau^*(\Gamma)=\tau^*(\Gamma[K])$.
%\newline
%\rm (и, разумеется, аналогичные равенства останутся выполненными,
%если заменить $K$ на любое большее конечное множество $K'\supseteq K$).

\Proof
We can take as $K$ the union of a minimal
cover and the support of a minimal fractional cover (which is finite
by the finite-fractional-cover lemma).

\smallskip

For example, we obtain the following fact.

\proclaim Infinite analogue of Aharoni--Holzman--Krivelevich's bound.
Suppose that the vertex set of a finitely covered finite-rank
hypergraph $\Gamma$
is partitioned into parts:
$V(\Gamma)=V_1\sqcup\dots\sqcup V_k$,
and $p_i=\max\limits_{e\in E(\Gamma)}|e\cap V_i|$.
Then
\item{\rm a)}
$
\tau(\Gamma)/\tau^*(\Gamma)
\le
\max\(p_1,\dots,p_k,{1\over2}\sum\limits_{i=1}^kp_i\);
$
moreover, the inequality is strict if $\sum\limits_{i=1}^kp_i\ge3$\;
\item{\rm b)}
if any $V_i$ contains no edge, then
$
\tau(\Gamma)/\tau^*(\Gamma)
\cases{
<
\rk\Gamma-1&for $k<\rk\Gamma$\;
\cr
\le
\(1-{1\over k}\)\rk\Gamma&for $k\ge\rk\Gamma$\;
\cr
}
$
\item{\rm c)}
if all $p_i$ equals 1, then
$
\tau(\Gamma)/\tau^*(\Gamma)
\le
\cases{
\rk\Gamma\cdot(k-\rk\Gamma+1)/k&for $k\ge(\rk\Gamma-1)\rk\Gamma$\;
\cr
\rk\Gamma\cdot k/(k+\rk\Gamma)+3-2\sqrt2&for $k<(\rk\Gamma-1)\rk\Gamma$.
}
$
\enditem

%\rem{
%Verify this. Ya them reformulated, simplifying, but could that-then
%наврать, especially|particularly in \`` boundary\"\ cases\dots
%}

\Proof
For finite hypergraphs, this strengthening of Lov\'asz's estimate [L75]
was proven in
[AHK96]:
\item{a)}
is Theorems 1 and 2;
\item{b)}
is Theorems 3 and 4;
\item{c)}
is Theorem 5 and simplifies (weakened) Theorem~6.
\enditem
It remains to apply the finite-trace theorem.

\smallskip

%%%%%%%%%%%%%%%%%%%%%%%%%%%%%%%
\s 3.
The cost of symmetry

A \emph{$G$-hypergraphs} is a hypergraphs equipped with an
action of a group $G$ by automorphisms.
The minimal size of a $G$-invariant fractional vertex cover of a
$G$-hypergraph $\Gamma$ (i.e., a fractional vertex cover $t$ such
that $t(g\o x)=t(x)$ for any $g\in G$ and $x\in X$) is denoted as
$\tau_G^*(\Gamma)$.

\Th 1.
For any finitely covered finite-rank $G$-hypergraph $\Gamma$
the equality
$\tau_G^*(\Gamma)=\tau^*(\Gamma)$ holds.
Moreover, for any finite-size fractional cover $t$,
the following functions are fractional covers too:
\item{\rm1)}
some finite cut $t[K]$
of $t$, where
$K\subseteq V$ is a finite set containing no point
with an infinite $G$-orbit;
\item{\rm2)}
the \emph{mean} $t_G\:v\mapsto{1\over|G\o v|}\sum\limits_{w\in G\o v}t(w)$
of $t$
\(where $1\over\hbox{\small an infinite cardinal}$ is assumed to be zero\).

\Proof
The equality $\tau_G^*(\Gamma)=\tau^*(\Gamma)$ follows immediately
from Assertion~2); indeed, the inequality $\ge$ is obvious, while the
inequality $\le$ holds, because the mean of $t$ is invariant and has size
at most $|t|$.

\goodbreak

\item{\rm1)}
\itemitem{--}
We can assume that the support $S$ of $t$ is finite
(because we can replace $t$ by its finite cut
the finite-fractional-cover lemma).

\itemitem{--}
We can assume that group $G$ leaves fixed all points of
$S$ with finite orbits
(because the finiteness of an orbit is equivalent to the finiteness
of the index of the stabiliser; so, we can replace $G$ with its
finite-index
subgroup~$\bigcap\bigl\{\St(s)\bigm| s\in S,\;|G:\St(s)|<\infty\bigr\}$,
and the points with infinite orbits remain points with infinite orbits).

\itemitem{--}
For each edge $e\in E(\Gamma)$,
there exists $g_e\in G$ such that $S\cap g_e\o e$
contains no point with
infinite orbit
(because otherwise
$
G=
\bigcup\limits_{s\in S, v\in e\atop |G\o v|=\infty}
\{g\in G\mid g\o v=s\}
$,
which is a finite union of cosets of stabilisers of
vertices $v$ with infinite orbits, i.e., the entire group $G$
is
covered by finitely many cosets of infinite-index
subgroups; and it is impossible
B.~Neumann's theorem [Neu54]).

\itemitem{--}
Now, for each edge $e\in E(\Gamma)$, we have
$
1\le\sum\limits_{v\in g_e\o e}t(v)
=
\sum\limits_{v\in e\cap S\atop|G\o v|<\infty}t(v)
$,
i.e., the cut $t[\{points\ with\ finite\ orbits\}]$
is a fractional cover, as required.

\item{\rm2)}
\itemitem{--}
We can assume that $t$ has a finite support $S$ containing no points with
infinite orbits
(by assertion 1), because $t_G(v)\ge\bigl(t[K]\bigr)_G(v)$
for all $v$).

\itemitem{--}
We can assume that the hypergraph $\Gamma$ is finite
(because $\Gamma$ can be replaced with its finite trace
$\Gamma\left[\bigcup\limits_{g\in G}g\o S\right]$;
on such a trace, the group $G$ acts, and
a fractional cover of such a trace gives a fractional cover of
$\Gamma$).

\itemitem{--}
Now, we can assume that the group $G$ is finite too
(because $G$ can be replaced with its image in the automorphism group
of $\Gamma$, which is a finite hypergraph already).

\itemitem{--}
When $G$ is finite,
the formula for the mean can be rewriten as:
$t_G(v)={1\over|G|}\sum\limits_{g\in G}t(g\o v)$.

\itemitem{--}
Now, everything is easy.
$\sum\limits_{v\in e}t(v)\ge1$ for each edge $e$.
Hence,
$1\le\sum\limits_{v\in g\o e}t(v)=\sum\limits_{v\in e}t(g\o v)$.
Summing over all~$g\in G$, we obtain
$
|G|\le\sum\limits_{g\in G}\sum\limits_{v\in e}t(g\o v)
=
\sum\limits_{v\in e}\sum\limits_{g\in G}t(g\o v)
$
and
$
1\le
\sum\limits_{v\in e}\sum\limits_{g\in G}{1\over|G|}t(g\o v)=
\sum\limits_{v\in e}t_G(v)
$
as required.
This completes the proof.

%\smallskip

\noindent
We call a $G$-hypergraph \emph{orbitally closed}
if, for any edge $e=\{v_1,\dots,v_l\}$
and any elements~$g_i\in G$
such that all vertices~$g_i\o v_i$ are different,
the set $\{g_1\o v_1,\dots,g_l\o v_l\}$
is an edge.
The \emph{orbit closure} of a $G$-hypergraph $\Gamma$ is the
(orbitally closed)
$G$-hypergraph $\=\Gamma$ with the same vertices as $\Gamma$
%$V(\=\Gamma)=V(\Gamma)$
and
$$
E(\=\Gamma)=
\bigl\{\{g_1\o v_1,\dots,g_l\o v_l\}\bigm| g_i\in G,\;
\{v_1,\dots,v_l\}\in E(\Gamma),\;
|\{g_1\o v_1,\dots,g_l\o v_l\}|=l
\bigr\}.
$$

\proclaim Orbit closure lemma.
For any
finitely covered finite-rank
$G$-hypergraph $\Gamma$,
\item{\rm1)}
$\tau(\Gamma)\le\tau(\=\Gamma)$
and
$\tau_G(\Gamma)=\tau_G(\=\Gamma)
\le
\tau(\=\Gamma)\cdot d(\Gamma)$,
where
$
d(\Gamma)=
\max\limits_{e\in E(\Gamma)\atop v\in V(\Gamma)}
\({|G\o v|\over|(G\o v)\setminus e|+1}\)
$,
\item{\rm2)}
and simpler relations for fractional covers holds:
$\tau^*(\Gamma)=\tau_G^*(\Gamma)=\tau_G^*(\=\Gamma)=\tau^*(\=\Gamma)$.

\Proof
\item{\rm1)}
The first inequality and equality are almost obvious: any cover
of $\=\Gamma$ is a cover of~$\Gamma$, and $G$-invariant covers
of $\Gamma$ and $\=\Gamma$ are the same things.
To prove the second
inequality, note that,
for a minimal cover $X$ of $\=\Gamma$,
$$
|(G\o v)\setminus X|
\le
\max\limits_{e\in E(\Gamma)}|e\cap G\o v|-1
\qbox{\sl
for any orbit $G\o v$ intersecting $X$.
}
\eqno{(**)}
$$
Indeed,
suppose the contrary
(i.e., $|(G\o v)\setminus X|\ge|e\cap G\o v|$
for an orbit $G\o v$ intersecting $X$
and any edge $e$);
then,
for each edge
$e=\{v_1,\dots,v_l\}\in E(\=\Gamma)$,
where $\{v_1,\dots,v_k\}=e\cap G\o v$, there are
elements $g_i\in G$ such that
$g_1\o v_1,\dots,g_k\o v_k\in(G\o v)\setminus X$
and $g_i\o v_i$ are all different;
%\newline
then
$E(\=\Gamma)\ni\=e=\{g_1\o v_1,\dots,g_k\o v_k,v_{k+1},\dots,v_l\}$
and
$
e\cap X\supseteq
\underbrace{\=e\cap X}_{\ne\emptyset}
=\{e_{k+1},\dots,e_l\}\cap X
\subseteq X\setminus(G\o v)
$,
i.e., $X\setminus(G\o v)$ is a cover of
$\=\Gamma$, which contradicts the minimality of $X$
(and proves $(**)$).
\newline
The union $Y$ of all orbits intersecting $X$ is an invariant
cover. Hence,
$$
\eqalign{
{\tau_G(\=\Gamma)\over\tau(\=\Gamma)}
&\le
{|Y|\over|X|}
\le
\max_v{|G\o v|\over|X\cap G\o v|}
=
\max_v{|G\o v|\over|G\o v|-|(G\o v)\setminus X|}
\lee^{**}
\max_{v,e}{|G\o v|\over|G\o v|-|e\cap G\o v|+1}
=
\cr
&=
\max_{v,e}{|G\o v|\over|(G\o v)\setminus e|+1}
\qbox{(where $\max\limits_v$ is the maximum over all vertices,
whose orbits
intersect $X$).}
}
$$
%\rem{Ya rewrote proof section 1), check.}

\item{\rm2)}
The first and last equalities hold by Theorem 1.
The second equality follows from that
$G$-invariant fractional covers for
$\Gamma$ and $\=\Gamma$ are the same things.
This completes the proof.

\smallskip

\enditem
These simple facts make it possible to obtain inequalities between
$\tau$ and $\tau_G$ from known results about
$\tau$ and $\tau^*$.
For example, we have the following theorem.

\proclaim Symmetric analogue of Aharoni--Holzman--Krivelevich's
estimate.
Suppose that the vertex set of a finitely covered
finite-rank $G$-hypergraph $\Gamma$
is partitioned into $G$-invariant parts:
$V(\Gamma)=V_1\sqcup\dots\sqcup V_k$,
$p_i=\max\limits_{e\in E(\Gamma)}|e\cap V_i|$,
and
$
d(\Gamma)=
\max\limits_{e\in E(\Gamma)\atop v\in V(\Gamma)}
\({|G\o v|\over|(G\o v)\setminus(e\cap G\o v)|+1}\).
$
Then
\item{\rm a)}
$
\tau_G(\Gamma)/\tau^*(\Gamma)
\le
d(\Gamma)\cdot
\max\(p_1,\dots,p_k,{1\over2}\sum\limits_{i=1}^kp_i\),
$
and the inequality is strict if $\sum\limits_{i=1}^kp_i\ge3$\;
\item{\rm b)}
if any $V_i$ contains no edge, then
$
\tau_G(\Gamma)/\tau^*(\Gamma)
\cases{
<
d(\Gamma)\cdot
(\rk\Gamma-1)&for $k<\rk\Gamma$;
\cr
\le
d(\Gamma)\cdot
\(1-{1\over k}\)\rk\Gamma&for $k\ge\rk\Gamma$;
\cr
}
$
\item{\rm c)}
if all $p_i$ equal 1, then
$
\tau_G(\Gamma)/\tau^*(\Gamma)
\le
\cases{
\rk\Gamma\cdot(k-\rk\Gamma+1)/k&for $k\ge(\rk\Gamma-1)\rk\Gamma$;
\cr
\rk\Gamma\cdot k/(k+\rk\Gamma)+3-2\sqrt2&
for $k<(\rk\Gamma-1)\rk\Gamma$.
}
$
\enditem
And the similar estimates hold for
$\tau_G(\Gamma)/\tau(\Gamma)$.

\Proof
\def\AHK{\mathrel{<\atop\le}}
Bounds for $\tau_G(\Gamma)/\tau(\Gamma)$ follow from bounds
for $\tau_G(\Gamma)/\tau^*(\Gamma)$, because
$\tau^*\le\tau$ always.

Let $\tau(\Gamma)/\tau^*(\Gamma)\AHK c(\Gamma,V_1,\dots,V_k)$ be
one of strict or nonstrict inequalities from the infinite analogue
of Aharoni--Holzman--Krivelevich's estimate
(i.e.,
$c(\Gamma,V_1,\dots,V_k)=
\max\(p_1,\dots,p_k,{1\over2}\sum\limits_{i=1}^kp_i\)$
in item~a), e.g.).
Then
$$
\tau_G(\Gamma)
\=^{\rm oc}
\tau_G(\=\Gamma)
\lee^{\rm oc}
d(\Gamma)\tau(\=\Gamma)
\AHK
c(\=\Gamma,V_1,\dots,V_k)d(\Gamma)\tau^*(\=\Gamma)
\=^{\rm oc}
c(\=\Gamma,V_1,\dots,V_k)d(\Gamma)\tau^*(\Gamma),
$$
where (oc) means a reference to the orbit-closure lemma.
It remains to note that
$c(\=\Gamma,V_1,\dots,V_k)=c(\Gamma,V_1,\dots,V_k)$,
and the conditions of c) implies that $d(\Gamma)=1$.
This completes the proof.

\smallskip

%\rem{
%This well, but guess завышенным --- there вот these
%two factor arise, is obtained\dots
%}

As a special (simplest) case of this theorem, we obtain the symmetric
analogue of Lov\'asz's estimate stated in the introduction.

%\proclaim Симметрический аналог оценки Ловаса.
%Если на $r$-дольном конечно покрытом гиперграфе $\Gamma$ ранга $r$
%действует группа $G$, сохраняя доли, то
%$
%\tau_G(\Gamma)/\tau(\Gamma)\le\tau_G(\Gamma)/\tau^*(\Gamma)\le r/2
%$,
%причём последнее неравенство строгое при $r>2$.

%%%%%%%%%%%%%%%%%%%%%%%%%%%%%
\s 4.
Open questions

It is known [AHK96] that Lov\'asz's estimate
is sharp
(as the most of its generalisations
from [AHK96]).
%(кроме случая $k<(\rk\Gamma-1)\rk\Gamma$ в пункте в),
%смотрите параграф 2).
The sharpness of the symmetric analogues is an open question.
We do not know the answer even in the simplest case,
i.e., for symmetric analogue of Lov\'asz's estimate for tripartite
3-hypergraphs.

\Question.
Is it true that, for any $\epsilon>0$, there exists a tripartite
hypergraph $\Gamma$ of rank three and a group $G$ acting on $\Gamma$
by automorphisms preserving the parts such that
$\tau_G(\Gamma)/\tau(\Gamma)>3/2-\epsilon$\?

In the following example, $\tau_G(\Gamma)/\tau(\Gamma)=4/3$.
Consider the 6-cycle $\Delta$ with coloured vertices and edges,
Figure 1 (r, g, and b mean red, green, and blue).

\goodbreak
%\vskip1cm plus 1cm minus5mm
\bigskip
\centerline{\input 1.PIC}
\nobreak%
%\vskip5mm%
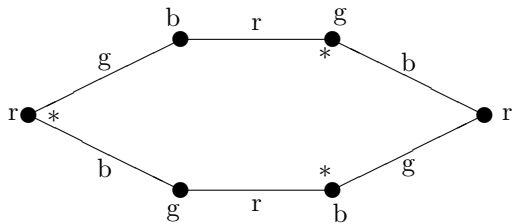
\centerline{Fig. 1}%
%\vskip1cm plus 1cm minus 5mm%
\goodbreak
\bigskip
%\par
%\noindent

Put $V(\Gamma)=V(\Delta)\cup E(\Delta)$ and
$E(\Gamma)=\bigl\{\{\hbox{an edge $e$ and two its endpoints}\}\mid
e\in E(\Delta)\bigr\}$. This 3-hypergraph $\Gamma$ is tripartite,
parts are denoted by colours.
The group $G$ (of order two) is the group of
colour-preserving automorphisms of $\Delta$ (the only
nontrivial automorphism is the central symmetry). Then $\tau(\Gamma)=3$ (a
minimal cover is the vertices marked by asterisk), and $\tau_G(\Gamma)=4$ (a
minimal invariant cover is the red and green vertices of $\Delta$).

\def\H#1{{\accent"7D #1}}
%%%%%%%%%%%%%%%%%%%%%%%
\References

[A85]
R. Aharoni,
Fractional matchings and covers in infinite hypergraphs,
Combinatorica, 5:3, (1985), 181-184.

[AH92]
R. Aharoni, R. Holzman,
Optimal fractional matchings and covers in infinite hypergraphs:
Existence and duality,
Graphs and Combinatorics, 8:1 (1992), 11-22.
%doi:10.1007/bf01271704

[AHK96]
R. Aharoni, R. Holzman, M. Krivelevich,
On a theorem of Lov\'asz on covers in $r$-partite hypergraphs,
Combinatorica, 16:2 (1996), 149--174.

[BZ19]
T. Banakh, D. van der Zypen,
Minimal covers of infinite hypergraphs,
Discrete Mathematics, 342:11 (2019), 3043-3046.
\arXiv 1808.08067

[DM58]
A. L. Dulmage, N. S. Mendelsohn,
Coverings of bipartite graphs,
Canadian Journal of Mathematics, 10 (1958), 517-534.

[F\"u88]
Z. F\"uredi,
Matchings and covers in hypergraphs,
Graphs and Combinatorics, 4:1 (1988), 115-206.

[KL21]
A. A. Klyachko, N. M. Luneva,
Invariant systems of representatives, or the cost of symmetry,
Discrete Mathematics, 344:6 (2021), 112361.
\arXiv 1908.03315

[KT25]
A. A. Klyachko, M. S. Terekhov,
Invariant systems of weighted representatives,
Journal of Algebraic Combinatorics, 61:3 (2025), 32.
\arXiv:2306.11883

[Ko19]
P. Komj\'ath,
Minimal covering sets for infinite set systems,
Discrete Mathematics, 342:6 (2019), 1849-1856.

[L75]
L. Lov\'asz,
A kombinatorika minimax t\'eteleir\H ol
(On the minimax theorems of combinatorics),
Matematikai Lapok 26:4 (1975), 209-264.

[L76]
L. Lov\'asz,
Covers, packings and some heuristic algorithms,
in: 5th British Comb.
Conf., Congressus Numerantium 15 Utilitas Math., Winnipeg, 1976, 417-429.

[Neu54]
B. H. Neumann,
Groups covered by permutable subsets,
J. London Math. Soc., s1-29:2 (1954), 236-248.

[T22]
M. S. Terekhov,
The cost of symmetry in connected graphs,
%Mat. Zametki, 112:6 (2022), 895-902.
Mathematical Notes, 112:6 (2022), 978-983.
\arXiv 2202.09590

[T25]
M. S. Terekhov,
The cost of symmetry for tailed stars,
arXiv:2510.03746~.

\end